\documentclass{article}
\usepackage{amssymb}
\usepackage{amsmath}

\usepackage{cite}

\usepackage{amssymb}
\usepackage{amsfonts}
\usepackage{amscd}
\usepackage[matrix]{xypic}

\newtheorem{theorem}{Theorem}

\newtheorem{corollary}[theorem]{Corollary}

\newtheorem{definition}[theorem]{Definition}

\newtheorem{lemma}[theorem]{Lemma}

\newtheorem{proposition}[theorem]{Proposition}

\begin{document}

\title{A class of nonassociative algebras}

\author{Michel Goze \thanks{%
corresponding author: e-mail: M.Goze@uha.fr} - Elisabeth
Remm \thanks{E.Remm@uha.fr.}\\
%EndAName
\\
{\small Universit\'{e} de Haute Alsace, F.S.T.}\\
{\small 4, rue des Fr\`{e}res Lumi\`{e}re - 68093 MULHOUSE - France}}
\date{}
\maketitle

A large class of nonassociative algebras, out of the Lie algebras has been studied, often with geometrical
structures and sometimes on fondamental spaces as the spaces of Hochschild cohomology. In this paper, we propose
to consider all these nonassociative algebras as algebras defined by the action of invariant subspaces of
the symmetric group $\Sigma_3$ on the associator of the considered laws.

\section{$\mathcal{V}$-algebras and $v$-algebras}

\subsection{Notations}

Let $\Sigma_3$ be the symmetric group of order $3$ and $\mathbb{K}\, [\Sigma_3]$ the associated group algebra,
where $\mathbb{K}$ is a field of characteristic zero.
 Let $\tau_{ij}$ be the transposition echanging $i$ and $j$, $c_1=(1,2,3)$ and $ c_2={c_1}^2$ the two $3$-cycles of $\Sigma_3.$
Every $v \in \mathbb{K}\,[\Sigma_3]$ is written
$$v=a_1 id + a_2 \tau _{12}  + a_3 \tau _{13} + a_4 \tau _{23}+a_5 c_1+a_6 c_2$$
or more generally
$$v=\sum_{\sigma \in \Sigma_3} a_{\sigma} \sigma
$$
where $a_{\sigma} \in \mathbb{K}.$

Consider the natural right action of $\Sigma_3$ on $\mathbb{K}\,[\Sigma_3]$
$$
\begin{array}{ccc}
\Sigma_3 \times \mathbb{K}\, [ \Sigma_3 ] & \rightarrow & \mathbb{K}\,[ \Sigma_3 ] \\
 (\sigma, \sum_i a_i \sigma_i ) & \mapsto & \sum_i a_i \sigma^{-1} \circ \sigma_i
\end{array}
$$
For every $v \in \mathbb{K}\,[\Sigma_3]$
we note $\mathcal{O}(v)$ the corresponding orbit of $v.$ Let $F_v=\mathbb{K}\,(\mathcal{O}(v))$ be the linear subspace
of $\mathbb{K}\,[\Sigma_3]$ generated by $\mathcal{O}(v).$ It is an invariant subspace of $\mathbb{K}\,[\Sigma_3].$
By Mashke's theorem it is a direct sum of irreducible invariant subspaces.

\subsection{$\mathcal{V}$-algebras and $v$-algebras}

Let $(\mathcal{A},\mu)$ be a $\mathbb{K}$-algebra where $\mu$ is the multiplication law of $\mathcal{A}.$
We note by $A_{\mu}$ the associator of $\mu$ that is
$$A_{\mu}= \mu \circ (\mu \otimes Id - Id \otimes \mu).
$$

Every $\sigma \in \Sigma_3$ defines a linear map noted $\Phi_{\sigma}$ given by
$$
\begin{array}{rccc}
\Phi_{\sigma} : & \mathcal{A}^{\otimes ^3} & \rightarrow & \mathcal{A}^{\otimes ^3} \\
& x_1 \otimes x_2 \otimes x_3 & \rightarrow & x_{\sigma^{-1}(1)} \otimes x_{\sigma^{-1}(2)} \otimes x_{\sigma^{-1}(3)}.
\end{array}
$$
If $v=\sum_{\sigma \in \Sigma_3} a_{\sigma} \sigma \in \mathbb{K}\,[\Sigma_3]$ we define the endomorphism
of $\mathcal{A}^{\otimes ^3}$ by putting:
$$\Phi_v=\sum_{\sigma \in \Sigma_3} a_{\sigma} \Phi_{\sigma}.$$

\begin{definition}
An algebra $(A,\mu)$ is a $\mathcal{V}$-algebra if there exists $v \in \mathbb{K}\, [\Sigma_3]$ such that
$$ A_{\mu} \circ \Phi_v =0.
$$
A $\mathcal{V}$-algebra is a $v$-algebra if for every $w \in \mathbb{K}\,[\Sigma_3]$ with $w \notin F_v$ we have that
$ A_{\mu} \circ \Phi_v =0$ and $ A_{\mu} \circ \Phi_w \neq 0$
\end{definition}
For a given vector $v$ in $\mathbb{K}\,[\Sigma_3]$ a vector $w \in F_v$ could exist such as $F_w$ is strictly contained in $F_v$. In this case a $v$-algebra is not a $w$-algebra. On the other hand the $w$-algbera is a $v$-algebra and for every vector $v' \in F_v$ such that $F_{v'}=F_v$, the notion of $v$ and $v'$-algebra are similar.

\medskip

\noindent Let us note that we have the following sequence for a $v$-algebra $(A,\mu)$
$$
\begin{CD}
\mathcal{A}^{\otimes ^3} @ >\Phi_v>> \mathcal{A}^{\otimes ^3} @ >A_{\mu}>> \mathcal{A}.
\end{CD}
$$
which represent a short complex.

\bigskip

\noindent {\bf Examples.}

1) Let $v=\sigma \in \Sigma_3.$ In this case $A_{\mu} \circ \Phi_v =0$ is equivalent to $A_{\mu}=0.$
This $v$-algebra is nothing other than an associative algebra.

2) Let $v=Id-\tau_{23}.$ A $v$-algebra satisfies
$$ A_{\mu}(x_1 \otimes x_2 \otimes x_3 -x_1 \otimes x_3 \otimes x_2)=0.
$$
Such an algebra is a Pre-Lie algebra \cite{G}.

\section{Irreducible $v$-algebras}

A $v$-algebra will be call irreducible if the $\Sigma_3$-module $F_v$ is irreducible. In this case we have
$dim \ F_v=1$ or $2$.

\subsection{One dimensional case}

\begin{proposition}
Let $v$ be in $\mathbb{K}\,[\Sigma_3]$ such that  $dim \ F_v=1$. Then $v$ is given, up to a multiplicative scalar, by one of the
two vectors $V$ and $W$ with
$$V= Id - \tau _{12} -  \tau _{13} - \tau _{23}+ c_1 +c_2,$$
$$W= Id + \tau _{12} +  \tau _{13} + \tau _{23}+ c_1 +c_2.$$
\end{proposition}
The first case correspond to the character of $\Sigma_3$ given by the signature, the second corresponding to the
trivial case.
We will keep the notation $V$ and $W$ to the end to designate these two found vectors. For a characterization of the
corresponding $v$-algebras we have need to recall the following definition.

\begin{definition}
A $\Bbb{K}$-algebra $\cal{A}$ with the multiplication $\mu$ is called a Lie-admissible if the bracket
$[x,y]=\mu(x,y)-\mu(y,x)$
satisfies the Jacobi identities, that is the algebra $(\cal{A},[,])$ is a Lie algebra.

\noindent A $\Bbb{K}$-algebra $(\cal{A},\mu)$  is called power-associative algebra if the associator $A_{\mu}$ satisfies
$$A_{\mu}(x,x,x)=0.$$
\end{definition}
The second condition $A_{\mu}(x,x,x)=0$ is equivalent to defining powers of a single element $x \in  \cal{A}$ recursively
by $x^1=x$, $x^{n+1}=\mu (x^n,x).$ An important class of power-associative algebra is the the class of Jordan algebras.

\begin{theorem}
There are two classes of irreducible $v$-algebras whose corresponding module $F_v$ is one-dimensional :

1. the Lie-admissible algebras given by $v=V$,

2. the power-associative algebras given by $v=W$.
\end{theorem}

\noindent {\it Proof.} Let $(\cal{A},\mu)$ be a $V$-algebra. Then we have $ A_{\mu} \circ \Phi_V =0. $ If we develop we
obtain:
$$ A_{\mu}(x,y,z)-A_{\mu}(y,x,z)-A_{\mu}(z,y,x)-A_{\mu}(x,z,y)+A_{\mu}(y,z,x)+A_{\mu}(z,x,y)=0.$$
This condition is equivalent to say that $[x,y]=\mu(x,y)-\mu(y,x)$ satisfies the Jacobi identity:
$$[x,[y,z]]+[y,[z,x]]+[z,[x,y]]=0.$$

\noindent Suppose now that $(\cal{A},\mu)$ is a $W$-algebra. Then, putting $[x,y]=\mu(x,y)-\mu(y,x)$ and
$\{x,y\}=\mu(x,y)+\mu(y,x)$, the condition $ A_{\mu} \circ \Phi_W =0  $ is equivalent to
$$[x,\{y,z\}]+[y,\{z,x\}]+[z,\{x,y\}]=0.$$
This last condition is just the linearization of the relation
$$A_{\mu}(x,x,x)=0.$$
In fact if $A_{\mu}(x,x,x)=0,$ we have also $A_{\mu}(x+y,x+y,x+y)=0$ which gives
$$A_{\mu}(x,x,y)+A_{\mu}(x,y,x)+A_{\mu}(y,x,x)+A_{\mu}(x,y,y)+A_{\mu}(y,x,y)+A_{\mu}(y,y,x)=0.$$
Moreover $A_{\mu}(x+y+z,x+y+z,x+y+z)=0$ implies taken account of the previous relations:
$$A_{\mu}(x,y,z)+A_{\mu}(x,z,y)+A_{\mu}(y,x,z)+A_{\mu}(y,z,x)+A_{\mu}(z,x,y)+A_{\mu}(z,y,x)=0$$
and the algebra is a $W$-algebra.
Then, if the characteristic of the field $\Bbb{K}$ is not $2$ (in fact we suppose it equal to $0$), the two last conditions
are equivalent.

\medskip

\noindent The classes of Lie-admissible algebras and power-associative algebras have been introduced by Albert in \cite{Al}.
\subsection{Two dimensional case, $\Bbb{K}=\Bbb{R}$}

In this section we suppose that the $\mathbb{K}\,[\Sigma_3]$-module $F_v$ is irreducible and $2$-dimensional. It is clear that, in this case,
$V$ nor $W$ do not belong to $F_v$.

\begin{proposition}
Let $v \in \mathbb{K}\,[\Sigma_3].$ If $F_v$ is a 2-dimensional irreducible invariant space  generated by
$v$ and $\tau_{12}(v)$ then
$$
\begin{array}{ll}
v &=  a_1Id+a_2\tau_{12}+( \alpha a_1 + \beta a_2) \tau_{13}-(\alpha a_1 +(1+\beta)a_2) \tau_{23}
\\
& + (\beta a_1 +\alpha a_2) c_1- ((1+\beta) a_1 +\alpha a_2) c_2,
\end{array}
$$
with $\alpha^2 =1+\beta+\beta ^2.$
\end{proposition}

\noindent {\it Proof.} Let $v=a_1Id +a_2 \tau _{12} +a_3  \tau _{13} +a_4 \tau _{23}+a_5 c_1 +a_6c_2$ be
in $\mathbb{K}\,[\Sigma_3]$ and suppose that $F_v$ is generated by $v$ and $\tau_{12}(v).$
As $\tau_{ij}(v) \in F_v$ the vector $v$ has to satisfied
$$
\sum_{i=1}^6 a_i =0.
$$
As $\left\{ \tau_{1i}\right\}_{ \ i=2,3}$ generates $\Sigma_3$, the $\mathbb{K}\,[\Sigma_3]$-module $F_v$ is $2$-dimensional if and only if $\{ v, \tau_{12}(v) \}$ are
independent and $\{  v,\tau_{12}(v), \tau_{13}(v) \}$ are related. Let $\tau_{13}(v)=\alpha v +\beta \tau_{12}(v)$ (let us recall that $\sigma (\sum a_i \sigma _i )= \sum a_i \sigma ^{-1} \circ \sigma _i )$ . This is
equivalent to the following system

$$
\left\{
\begin{array}{l}
Aa_1=Ba_2 \ ; \  Aa_3=Ba_6 \ ; \  Aa_4=Ba_5 \\
Aa_2=Ba_1 \ ;  \ Aa_6=Ba_3 \ ;  \ Aa_5=Ba_4 \\
A=1-\alpha^2+\alpha^2\beta-\beta^3 \\
B=\alpha -\alpha^3+\alpha \beta^2+\alpha \beta
\end{array}
\right.
$$
which implies that $(A^2-B^2)a_1=(A^2-B^2)a_2=0.$

\medskip

\noindent $1^{{\rm srt}}$ case. $A^2-B^2 \neq 0.$ Then we have $a_1=a_2=0.$ As $\tau_{ij}(v) \in F_v$ for $(i,j) \neq (1,2),$
we deduce  $v=0$ which is impossible.

\smallskip

\noindent $2^{{\rm nd}}$ case. $A^2-B^2=0$.

\smallskip

i) $A=B$ that is $1-\alpha^2+\alpha^2\beta -\beta^3=\alpha(1-\alpha^2+\beta^2+\beta)$. We deduce
$A(a_1-a_2)=0.$ Then if $A \neq 0,$ we obtain  $a_1=a_2,a_3=a_6,a_4=a_5$ and in this case $v=\tau_{12}(v)$
which contradicts the hypothesis. Thus $A=B=0$ and the coefficients $\alpha$ and $\beta$ satisfy
$$
\left\{
\begin{array}{l}
1-\alpha^2+\alpha^2+\alpha^2\beta -\beta^3=0 \\
\alpha(1-\alpha^2+\beta^2+\beta)=0
\end{array}
\right.
$$
If $\alpha=0, \ \beta=0$ (here we choose $\mathbb{K}=\mathbb{R}$), we have $W \in F_v$ and $F_v$ is not irreducible.
Thus $\alpha \neq 0$ and
$$
\left\{
\begin{array}{l}
1-\alpha^2=-\beta -\beta^2 \\
1-\alpha^2=-\beta(\alpha^2-\beta^2)
\end{array}
\right.
$$
If $\beta=0, \ \alpha=-1$ or $\alpha=1.$

\begin{lemma}
If $\beta=0$ and $\alpha=1$ then $v=(a_1,a_2,a_1,-a_1-a_2,a_2,-a_1-a_2).$
\end{lemma}

\noindent {\it Proof.}
We write in this case $\tau _{13}(v)=v$ which implies
$$v=(a_1,a_2,a_1,a_4,a_2,a_4)
$$
As $\tau _{23}(v)=\alpha'v+\beta'\tau _{12}(v),$ we deduce that
$$(\alpha'-\beta')(a_2-a_1)=0.
$$
If $a_1=a_2,$ $v=(a_1,a_1,a_1,a_4,a_1,a_4)$ with
$$
\left\{
\begin{array}{l}
a_4=(\alpha'+\beta')a_1 \\
a_1=\alpha'a_1+\beta'a_4=\alpha'a_4 +\beta'a_1
\end{array}
\right.
$$
which gives that $\beta'=-1$ or $a_4=a_1.$ The second case corresponds to the vecteur $W.$ Thus we have
$$
v=(a_1,a_1,a_1,-2a_1,a_1,-2a_1)
$$
because $-2-\alpha'+\alpha'^2=0$ which means that $\alpha'=-1$ or $\alpha'=2$ but for this last choice $v=W.$

\smallskip

\noindent If $a_1 \neq a_2$ then $\alpha'=\beta'$ and we deduce that
$$
\left\{
\begin{array}{l}
a_4=\alpha'(a_1+a_2) \\
a_2=\alpha'(a_1+a_4) \\
a_1=\alpha'(a_2+a_4)
\end{array}
\right.
$$
which implies that $(a_1+a_2+a_4)(1-2\alpha')=0$ and $a_4=-a_1-a_2.$ The vector $v$ has the expecting form. If $\alpha'=\frac{1}{2}$ then $v=W.$
Whence the announced lemma.

\begin{lemma}
If $\beta=0$ and $\alpha=-1$ then $v=(a_1, a_2,-a_1, a_1-a_2, -a_2, -a_1+a_2).$
\end{lemma}

\noindent {\it Proof.} The calcul is similar to the previous proof.

\medskip

\noindent Suppose now that $\beta \neq 0$. In this case the system between $\alpha$ and $\beta$ is equivalent to the single equation
$$\alpha ^2 = 1 + \beta + \beta ^2.$$
We deduce $a_3= \alpha a_1 + \beta a_2$, $ a_5 = \alpha a_2 + \beta a_1$, and $\beta a_6 = (1-\alpha ^2)a_1-\alpha \beta a_2.$
As $\alpha ^2 = 1 +\beta + \beta ^2$ and $\beta \neq 0$, we deduce $a_6 = -(1+ \beta )a_1 - \alpha a_2$. Likewise we have
$a_4=-\alpha a_1 - (1+\beta ) a_2$. The vector $v$ has the following form :
$$v=(a_1,a_2,\alpha a_1+\beta a_2, -\alpha a_1-(1+\beta ) a_2, \beta a_1 +\alpha a_2 ,-(1+\beta ) a_1-\alpha a_2 ),$$
with $\alpha ^2 = 1 + \beta + \beta ^2.$

\medskip

ii) $A=-B$. If $A \neq 0$, then $a_1=-a_2$, $a_3=-a_6$, and $a_4=-a_5.$ This implies $v= \tau _{12}(v),$ which contradicts the hypothesis. Then
$A=0$ and we come back in the previous case. $\clubsuit$

\bigskip

\noindent {\bf Consequence.} Let $(A,\mu)$ be a $v$-algebra with $dim \ F_v=2$ and $F_v$ irreducible. We have one of the
following cases
$$ F_v=\mathbb{K}\left\{ v,\tau(v) \right\} \ {\rm or } \ F_v=\mathbb{K}\left\{ v,c(v) \right\} $$
where $\tau$ is a transposition and $c$ a 3-cycle. But all the transpositions are conjugated in $\Sigma_3$
then if $F_v= \mathbb{K}\left\{ v,\tau(v) \right\}$ there exists $v' \in F_v$ such that
$F_{v'}=F_{v}=\mathbb{K}\left\{ v',\tau_{12}(v') \right\}$. In this case the expression of $v'$ is given by
the proposition 5. Now the second case $F_v=\mathbb{K}\left\{ v,c(v) \right\}$ can be discussed. In fact
as the transpositions generate $\Sigma_3$, the dependance of the vectors $\left\{ v,\tau(v) \right\}$ for every
tranposition implies the dependance of the vectors $\left\{ v,c(v) \right\}$ for every $3$-cycle.
Then every $v$-algebra with $dim \ F_v=2$ and $F_v$ irreducible is discribed by the vector $v'$
such that $F_{v'}=\mathbb{K}\left\{ v',\tau_{12}(v') \right\}.$ Thus we can take $v$ given by the
propositon $5.$ We have
$$
\begin{array}{ll}
v &=  a_1Id+a_2\tau_{12}+( \alpha a_1 + \beta a_2) \tau_{13}-(\alpha a_1 +(1+\beta)a_2) \tau_{23}
\\
& + (\beta a_1 +\alpha a_2) c_1- ((1+\beta) a_1 +\alpha a_2) c_2 ,
\end{array}
$$
with $\alpha^2 =1+\beta+\beta ^2,$ and
$$
\begin{array}{ll}
\tau _{12}(v) &=  a_2Id+a_1\tau_{12}- ((1+\beta) a_1 +\alpha a_2) \tau_{13}+ (\beta a_1
+\alpha a_2) \tau_{23}
\\
& -(\alpha a_1 +(1+\beta)a_2) c_1+( \alpha a_1 + \beta a_2) c_2 .
\end{array}
$$
Let us consider the vector $u=v-\tau_{12}(v).$ We have
$$
\begin{array}{ll}
\tau_{12}(u) & =-u, \\
\tau_{13}(u) & =\tau_{13}(v)-c_1(v) \\
 & =\alpha v +\beta \tau_{12}(v) -\beta(v)-\alpha \tau_{12}v=(\alpha - \beta)(u) \\
\tau_{23}(u) & =\tau_{23}(v)-c_2(v)=-\alpha v-(1+\beta)\tau_{12}(v)-(1+\beta) v-\alpha\tau_{12}(v)\\
& =-(\alpha+\beta+1)(v+\tau_{12}(v))
\end{array}
$$
Thus $\left\{ u,\tau_{23}(u) \right\}$ generate $F_v$ that is $F_v=F_u$ and every $v$-algebra is an $u$-algebra.
But the vector $u$ as the simplified following form
$$
\begin{array}{ll}
u & =  (a_1-a_2,a_2-a_1,(\alpha + \beta +1)a_1+(\alpha + \beta )a_2,-(\alpha + \beta )a_1-
(\alpha + \beta +1)a_2, \\
& (\alpha + \beta )a_1+(\alpha + \beta +1)a_2,-(\alpha + \beta +1)a_1-(\alpha + \beta )a_2) \\
& = (\lambda_1,-\lambda_1,\lambda_2,-\lambda_3,\lambda_3,-\lambda_2)
\end{array}
$$
with $\lambda_1+\lambda_3-\lambda_2=0,$ that is
$$u=(\lambda_1,-\lambda_1, \lambda_1+\lambda_3,-\lambda_3,\lambda_3,-\lambda_1-\lambda_3).
$$

\begin{theorem}
Every $v$-algebra with $F_v$ irreducible and 2-dimensional is given by the following identity:
$$
\begin{array}{l}
\lambda_1 (A_{\mu}(x,y,z)- A_{\mu}(y,x,z))+(\lambda_1 + \lambda_3)(A_{\mu}(z,y,x) -A_{\mu}(z,x,y) \\
-\lambda_3(A_{\mu}(x,z,y)-A_{\mu}(y,z,x))=0
\end{array}
$$
\end{theorem}

\section{Lie-admissible $\mathcal{V}$-algebras}

The aim of this section is to describe all the $v$-algebras which are Lie-admissible.

\begin{lemma}
A $v$-algebra is Lie-admissible if and only if $V \in F_v.$
\end{lemma}

\noindent {\it Proof.}

\smallskip

\noindent $(\Rightarrow)$ We have that
$$
\begin{array}{lll}
A_{\mu} \circ \Phi_v=0 & \Rightarrow & A_{\mu} \circ \Phi_{\sigma(v)}=0 \ {\rm for \ all } \ \sigma \in \Sigma_3 \\
& \Rightarrow &  A_{\mu} \circ \Phi_{v'}=0 \ {\rm for \ all } \ v' \in F_v.
\end{array}
$$
Thus we have that $V \in F_v \Rightarrow A_{\mu} \circ \Phi_V=0$ and the $v$-algebra is Lie-admissible.

\smallskip

\noindent $(\Leftarrow)$
By definition of the $v$-algebra $A_{\mu} \circ \Phi_{v'}=0$ if and only if $v' \in F_v$ Then if the $v$-algebra
is Lie-admissible the vector $V$ is in $F_v$. $\clubsuit$

\bigskip

Let us first recall some results on the particular Lie-admissible $v$-algebras which are the $G_i$-algebras.

\subsection{The $G_i$-algebras of \cite{G.R}}
In \cite{G.R} and \cite{R} the class of Lie-admissible algebras related with the subgroups of
$\Sigma _3$ has been classified. We recall briefly this result.
\begin{definition}
Let $G$ be a subgroup of $\Sigma _3$. We call $G$-associative algebra an algebra $(\cal{A},\mu)$ satisfying
$$\sum _{\sigma \in G}A_{\mu} \circ \Phi _{\sigma} = 0.$$
\end{definition}

We have denoted $G_1={Id}, G_2=\left\{Id,\tau_{12}\right\},G_3=\left\{Id,\tau_{23}\right\},G_4=\left\{Id,\tau_{13}\right\},$
$G_5=\left\{Id,c_1,c_2\right\},G_6=\Sigma_3$ the subgroups of $\Sigma_3.$
When $G=G_6$ the class of $\Sigma _3$-associative algebras is the full class of Lie-admissible algebras,
when $G=G_1$ the corresponding class is the full class of associative algebras. For other cases we obtain the class of Vinberg algebras for $G_2$, the pre-Lie algebras for $G_3$. Each one of these algebras can be
presented as a $v$-algebra :

- the $\Sigma _3$-associative algebras are the $V$-algebras.

- the  ($G_1$-)associative algebras are the $v$-algebras with $v=Id.$

- the Vinberg algebras, i.e. the $G_2$-associative algebras are the $v$-algebras corresponding to $v=Id - \tau _{12}.$

- the pre-Lie algebras, i.e. the $G_3$-associative algebras are the $v$-algebras corresponding to $v=Id - \tau _{23}.$

- the $G_4$-associative algebras are the $\{Id - \tau _{13}\}$-algebras.

- the $G_5$-associative algebras are the $\{Id + c_1 + c_2\}$-algebras.

\noindent We will generalize this list considering not only the invariant spaces $F_v$ associated to a subgroup
of $\Sigma _3$, but all the invariant spaces.

\subsection{Classification  of Lie-admissible $\mathcal{V}$-algebras}

\begin{theorem}
Every Lie-admissible $v$-algebra corresponds to the following one:

$\bullet$ type ($I$): $dim \ F_v=1$ and $\mu(x,y)=x.y$ satisfies:
$$
\begin{array}{l}
(x.y).z-x.(y.z)-(y.x).z+y.(x.z)-(x.z).y+x.(z.y)-(z.y).x+z.(y.x) \\
\qquad \qquad +(y.z).x-y.(z.x)+(z.x).y-z.(x.y)=0.
\end{array}
$$
This identity defines the category of Lie-admissible algebras.

\medskip

$\bullet$ type ($II$): $dim \ F_v=2$ and the $v$-algebra is a power-associative algebra given by
$$(x.y).z-x.(y.z)+(y.z).x-y.(z.x)+(z.x).y-z.(x.y)=0$$
(that is a $G_5$-associative algebra).

\medskip

$\bullet$ type ($III$): $dim \ F_v=3$. The $v$-algebra corresponds to
$$
\begin{array}{l} (III_1) \ (x.y).z-x.(y.z)+t[(y.x).z-y.(x.z)]-[(x.z).y-x.(z.y)] \\
\qquad \qquad -t[(z.x).y-z.(x.y)]=0,
\end{array}
$$
where $t \neq 1 $ or by
$$
\begin{array}{l}
(III_2) \ (x.y).z-x.(y.z)-(z.y).x+z.(y.x). \qquad  \qquad \qquad \qquad \quad
\end{array}
$$
or by
$$
\begin{array}{l}
(III_3) \ (x.y).z-x.(y.z)-(y.x).z+y.(x.z)-2(x.z).y+2x.(z.y) \\
\qquad \qquad +2(y.z).x-2y.(z.x)=0.
\end{array}
$$

\medskip

$\bullet$ type ($IV$): $dim \ F_v=4.$ The $v$-algebras are of the following type.
$$
\begin{array}{l}
(IV_1) \ 2[(x.y).z-x.(y.z)]+(1+t)[(y.x).z-y.(x.z)]+(z.y).x \\
\qquad -z.(y.x)+[(y.z).x-y.(z.x)]+(1-t)[(z.x).y-z.(x.y)]=0
\end{array}
$$
with $t \neq 1$,
$$
\begin{array}{l}
(IV_2) \ 2[(x.y).z-x.(y.z)]+[(y.x).z-y.(x.z)]+[(x.z).y-x.(z.y)] \\
\qquad \qquad +[(y.z).x-y.(z.x)]+[(z.x).y-z.(x.y)]=0
\end{array}
$$
$$
\begin{array}{l}
(IV_3) \ 2[(x.y).z-x.(y.z)]+[(z.y).x-z.(y.x)]-[(x.z).y-x.(z.y)] \\
\qquad \qquad +3[(y.z).x-y.(z.x)]+[(z.x).y-z.(x.y)]=0
\end{array}
$$

\medskip

$\bullet$ type ($V$): $dim \ F_v=5$ and we have
$$
\begin{array}{l}
2[(x.y).z-x.(y.z)]-(y.x).z+y.(x.z)-(z.y).x+z.(y.x) \\
\qquad \qquad -((x.z).y+x.(z.y))+((y.z).x-y.(z.x))=0
\end{array}
$$

\medskip

$\bullet$ (type $VI$): $dim \ F_v=6.$ This correspond to the class of associative algebras
$$(x.y).z-x.(y.z)=0
$$

\end{theorem}

\noindent {\it Proof.}

\medskip

\noindent 1) $dim \ F_v=1.$ This corresponds to theorem 4.1.

\medskip

\noindent 2) $dim \ F_v=2$. As $V \in F_v,$ $F_v=F_V \oplus F_W.$ We can take $v=\frac{1}{2}(V+W)$ that is
$v=(1,0,0,0,1,1).$ Then $\mu$ satisfies
$$A_{\mu}(x,y,z)+A_{\mu}(y,z,x)+A_{\mu}(z,x,y)=0.$$

\medskip

\noindent 3) $dim \ F_v=3$. In this case $F_v=F_V \oplus F_u$ with $dim \ F_u=2$ and $F_u$ irreducible because we have only two non equivalent irreducible 1-dimensional representations. We have seen that $u$ can be chosen of the form
$$u=(\lambda _1,-\lambda _1,\lambda_1 + \lambda_3, -\lambda _3, \lambda _3, -\lambda_1 -\lambda_3).$$
Suppose that $\lambda_1 \neq 0$ and consider the vector $v'=u-\lambda_1 V:$
$$v'=(0,0,2\lambda_1+\lambda_3, \lambda_1-\lambda_3, -\lambda_1+\lambda_3,-2\lambda_1-\lambda_3).$$
Then $\mathbb{K}\,\mathcal{O}(v')=\mathbb{K}\left\{ v',\tau_{13}(v'),\tau_{23}(v') \right\}.$ This space is
of dimension 3 if $(b \neq a)$ or $(a=-b \neq 0)$ with $a=2\lambda_1+\lambda_3$ and $b=-\lambda_1+\lambda_3.$
The first case correspond
to $\lambda_1 \neq 0$ and the second case to $\lambda_1=0, \lambda_3 \neq 0.$ If $\lambda_1 \neq 0,$ we can write
$$v'=(0,0,a,-b,b,-a)$$ with $a \neq b.$ To simplify, using $v"=\tau_{13}(v'),$ we have
$\mathbb{K} \ \mathcal{O}(v")=\mathbb{K} \ \mathcal{O}(v')=F_v$ and $v"=(a,b,0,-a,0,-b).$ If $a \neq 0,$ divising by $a$
and putting $t=\frac{b}{a},$ thus we can come down to $v=(1,t,0,-1,0,-t)$ with $t \neq 1$
($t=1$ corresponds to $\lambda_1=0$).
This corresponds to the case $(III_1).$ If $a=0$ and $b \neq 0,$ we have $v=(0,-1,0,0,0,1)$
and we obtain the case $(III_2).$ If $\lambda_1=0,$ we get
$u=(0,0,\lambda_3,-\lambda_3,\lambda_3,-\lambda_3)=\lambda_3(0,0,1,-1,1,-1).$ We take $u+V=(1,-1,0,-2,2,0)$ and
we find the case $(III_3).$

\medskip

\noindent 4) $dim \ F_v=4.$ As $F_V \subset F_v,$ $F_v=F_u \oplus F_V \oplus F_W$ with $dim \ F_u=2$
and $F_u$ irreducible. The subspace $F_u \oplus F_V$ corresponds to the previous case.
Let $v=(1,t,0,-1,0,-t)$ a vector associated to an algebra of type $(III_1).$
In these conditions $t \neq 1.$ Let $\widetilde{v}=v+W=(2,1+t,1,0,1,1-t).$ Then $dim \ F_{\widetilde{v}}=4$
(because $t \neq 1$) and
$W=\frac{1}{2}(\tau_{23}(\widetilde{v})+\widetilde{v}), V=(1-t)(-\frac{1}{2}\widetilde{v}-\frac{3}{2}\tau_{12}(\widetilde{v})+\tau_{13}(\widetilde{v})+\tau_{23}(\widetilde{v})).$
Thus the $v$-algebras correspond to the equations
$$
\begin{array}{rl}
2 A_{\mu}(x,y,z)+(1+t)A_{\mu}(y,x,z)+A_{\mu}(z,y,x) & \\
+A_{\mu}(y,z,x)+(1-t)A_{\mu}(z,x,y) & =0
\end{array}
$$
These algebras are Lie-admissible and power-associative.

If $v$ is generating an algebra of type $(III_2)$ then
$$v=(1,0,-1,0,0,0)$$
In this case $\widetilde{v}=v+W=(2,1,0,1,1,1).$

Finally if $v=(1,-1,0,-2,2,0),$ the vector generating the algebras of type $(IV_3)$ is of type
$\widetilde{v}=(2,0,1,-1,3,1).$

\medskip

\noindent 5) $dim \ F_v=5.$ In this case $F_v=F_{v_1}\oplus F_{v_2} \oplus F_V$ with $F_{v_1}$ and $F_{v_2}$
irreducible and $2$-dimensional. As $dim \ F_v=5,$ the vectors
$v, \tau_{12}(v),\tau_{13}(v),\tau_{23}(v),c_1(v),c_2(v)$ are of rank 5. There exists a linear relation between
these vectors. Putting
$v=\alpha_1 Id + \alpha_2 \tau_{12}+\alpha_3 \tau_{13}+\alpha_4 \tau_{23}+\alpha_5 c_1 +\alpha_6 c_2,$ we have obviously
$$v +\tau_{12}(v)+\tau_{13}(v)+\tau_{23}(v)+ c_1(v)+ c_2(v)=\sum_{i=1}^6 \alpha_i W.$$
But $W \notin F_v$ then $\sum_{i=1}^6 \alpha_i=0.$ We find again the characterization of invariant subspace of $\mathbb{K}\,[\Sigma_3]$ of codimension one:

 "If $dim \ F_v=5$ then $v$ satisfies $\sum \alpha_i =0."$

\noindent Then there exists only one type of such subspace. A basis of $F_v$ can be obtained by the vectors
$\left\{ e_1=(1,0,0,0,0,-1), e_2=(0,1,0,0,0,-1), e_3=(0,0,1,0,0,-1),\right.$

\noindent $\left. e_4=(0,0,0,1,0,-1), e_5=(0,0,0,0,1,1) \right\}.$

Consider for example the vector $e_1$ and let us compute $\mathcal{O}(e_1).$ We obtain
$$\mathcal{O}(e_1)=\left\{e_1,\tau_{12}(e_1),\tau_{13}(e_1),c_1 (e_1) \right\}$$
and $dim \ \mathbb{K} \,( \mathcal{O}(e_1))=4.$ But $V \notin \mathbb{K} \,( \mathcal{O}(e_1))$ and
$dim \ \mathbb{K} \,( \mathcal{O}(e_1+V))=5.$ The unicity of $F_v$ implies that we can reduce the vector $v$ to $e_1+V.$
Then the class of Lie-admissible $v$-algebras with $dim \ F_v=5$ is given by $v=(2,-1,-1,-1,1,0)$
that is it satisfies
$$
\begin{array}{ll}
(V): \ &  2x.(y.z) -2(x.y).z-y.(x.z)+(y.x).z-z.(y.x)+(z.y).x \\
& -x.(z.y)+(x.z).y+y.(z.x)-(y.z).x=0.
\end{array}
$$

\medskip

\noindent 6) $dim \ F_v=6:$ If $dim \ F_v=6$ then $F_v=\mathbb{K}\,[\Sigma_3].$ Let us consider the vector
$e_1=(1,0,0,0,0,0).$ We have $\mathbb{K} \,( \mathcal{O}(e_1))=\mathbb{K}\,[\Sigma_3]=F_v.$
Then $v$ can be reduced to the vector $e_1.$ In this case a Lie-admissible $v$-algebra satisfies
$x.(y.z)-(x.y).z=0$ that is, we find again the class of associative algebras.

\bigskip

\noindent{\it Remark. }The correspondance between the $v$-algebras and the $G_i$-associative algebras is the following:
 type $(I)$ corresponds to $G_6$, type $(II)$ to $G_5$, type $(III_1)$ to $G_2$ when $t=-1$ and to $G_3$ when $t=0$, type $(III_2)$
 corresponds to $G_4.$

\section{Power-associative $v$-algebras}

\subsection{Classification}

We have seen that the category of power-associative algebras corresponds to the class of the $W$-algebras.
In the previous section some classes of Lie-admissible algebras are also power-associative algebras.
These cases correspond to $v$-algebras which $F_v$ contains $F_V \oplus F_W.$ This appears if

\noindent $\bullet \ dim \ F_v=2: \ F_v=F_V \oplus F_W;$

\noindent $\bullet \ dim \ F_v=4: \ F_v=F_u \oplus F_V \oplus F_W;$

\noindent $\bullet \ dim \ F_v=6: \ F_v=\mathbb{K}\,[\Sigma_3].$

\medskip

Then it remains only the cases $dim \ F_v=3$ and $dim \ F_v=5.$

\smallskip

\noindent {\it First case.} We have $dim \ F_v =3$ and $F_v=F_u \oplus F_W$ with $F_u$ irreducible and $2$-dimensional.
The vector $u$ is of the form
$u=(\lambda_1,-\lambda_1, \lambda_1+\lambda_3, -\lambda_3, \lambda_3, -\lambda_1-\lambda_3).$
If $\lambda_1 \neq 0,$ then $v=u-\lambda_1 W$ satisfies $dim \ \mathbb{K} \,( \mathcal{O} \ (v))=3$ as soon as
$\lambda_1 \neq 0$ or $\lambda_1+\lambda_3 \neq 0.$ If these conditions are satisfied,
$\widetilde{v}=(0,-2,t,-1-t,-1+t,-2-t)$ generates $\mathbb{K} \,( \mathcal{O} \, (v)).$

\noindent If $\lambda_1=0$ then $u=(0,0,\lambda_3,-\lambda_3,\lambda_3,-\lambda_3)$ and
$\widetilde{v}=\frac{u}{\lambda_3}+ W=(1,1,1,0,1,0)$ generates $F_v.$

\smallskip

\noindent {\it Second case.} $dim \ F_v =5.$ Then $F_v=F_{u_1} \oplus F_{u_2} \oplus F_W$ with $F_{u_i}$ irreducible
and $2$-dimensional for $i=1,2.$ If $v=(a_1,a_2,a_3,a_4,a_5,a_6)$ then
$$
\begin{array}{l}
v-\tau_{12}(v)-\tau_{13}(v)-\tau_{23}(v)+c_1(v)+c_2(v)=(a_1-a_2-a_3-a_4+a_5+a_6)V.
\end{array}
$$
As $V\notin F_v,$ we have that $a_1-a_2-a_3-a_4+a_5+a_6=0$ and
$$ F_v=\left\{  (a_1,a_2,a_3,a_4,a_5,a_6) \ / a_1-a_2-a_3-a_4+a_5+a_6=0 \right\}.$$
Let $v_1=(1,0,0,0,0,-1)$ be in $F_v.$ We have proved in the last section that $\mathbb{K} \,( \mathcal{O}(v_1))$
is of dimension 4. Then $v=v_1 +W$ is a "generator" of $F_v.$ As $v=(2,1,1,1,1,0),$ te corresponding
power-associative $v$-algebra satisfies
$$
\begin{array}{l}
2 A_{\mu}(x,y,z)+A_{\mu}(y,x,z)+A_{\mu}(z,y,x)+A_{\mu}(x,z,y)+A_{\mu}(y,z,x)=0.
\end{array}
$$

\begin{theorem}
Every power-associative $v$-algebra corresponds to the following one:

$\bullet$ type $(I')$: $dim \ F_v=1$ and $\mu(x,y)=x.y$ satisfies:
$$
\begin{array}{l}
(x.y).z-x.(y.z)+(y.x).z-y.(x.z)+(x.z).y-x.(z.y)+(z.y).x-z.(y.x) \\
\qquad \qquad +(y.z).x-y.(z.x)+(z.x).y-z.(x.y)=0.
\end{array}
$$
This identity defines the category of power-associative algebras.

\medskip

$\bullet$ type $(II)$ $(dim \ F_v=2)$ and we have the algebras defined in theorem 11.

\medskip

$\bullet$ type $(III')$: $dim \ F_v=3$. The $v$-algebra correspond to

$(III'_{\ 1}) :$
$$
\begin{array}{l}  -2[(x.y).z-x.(y.z)]-(2+t)[(z.y).x-z.(y.x)]+(t-1)[(x.z).y-x.(z.y)] \\
\qquad \qquad -(1+t)[(y.z).x-y.(z.x)]+t[(z.x).y-z.(x.y)]=0,
\end{array}
$$
$(III'_{\ 2}) :$
$$
\begin{array}{l}
 (x.y).z-x.(y.z)+(y.x).z-y.(x.z)+(z.y).x-z.(y.x)+(y.z).x-y.(z.x).
\end{array}
$$

\medskip

$\bullet$ type $(IV)$: $dim \ F_v=4.$ The $v$-algebras are of type $(IV_1)$ or $(IV_1)$ of the theorem 11.

$\bullet$ type $(V')$: $dim \ F_v=5$ and we have
$$
\begin{array}{l}
2[(x.y).z-x.(y.z)]+(y.x).z-y.(x.z)+(z.y).x-z.(y.x) \\
\qquad \qquad +((x.z).y-x.(z.y))+((y.z).x-y.(z.x))=0
\end{array}
$$
 This class corresponds to the alternative algebras.
\medskip

$\bullet$ type $(VI)$: $dim \ F_v=6$ and we have the class of associative algebras.

\end{theorem}

Let us examine more paticularly type $(V')$.

\begin{definition} An algebra $(\mathcal{A}, \mu)$ is alternative if its product satisfies
$$
A_{\mu} (x,x,y)=A_{\mu} (y,x,x)=0
$$
\end{definition}

This condition is equivalent to the following system
$$
A_{\mu}\circ \Phi_{v_1}=A_{\mu}\circ \Phi_{v_2}=A_{\mu}\circ \Phi_{v_3}=0
$$
with $v_1=Id+\tau_{12},v_2=Id-c_1,v_3=Id+\tau_{13.}$ The vectors are in the orbit of the vector $v$ associated to the type $(V')$ and in fact $\mathbb{K}\,(\mathcal{O}\, \left\{ v_1,v_2,v_3\right\})=\mathcal{K}\, (\mathcal{O} \, (v)).$ So the calss of algebras of type $(V')$ corresponds to the class of alternative algebras.

\subsection{A characteristic example: the octonions algebra}

The octonions algebra also called Caley algebra is a $8$-dimensional algebra which is a $W$-algebra that is a  power-associative $\mathcal{W}$-algebra of type $(V')$. More generally the $8$-dimensional Caley-Dixon algebras of composition are $W$-algebras of type $V'$.

Recall that any division composition algebra over $\mathbb{R}$
or $\mathbb{C}$ or the quaternions's algebras or the octonions's algebras.

\section{Tensor products of $\mathcal{V}$-algebras}

It is wellknown that the category of associative algebras is stable by tensor product. We saw in \cite{G.R} that it is
not the case for the $G_i$-associative algebras which are not associative. Here we will show that

\begin{proposition}
Let $A$ and $B$ be two $\mathcal{V}$-algebras. Then $A \otimes_{\mathbb{K}}B$ is a $\mathcal{V}$-algebra
if and only $A$ and $B$ are associative.
\end{proposition}

\noindent {\it Proof.} Let $A_\mu$ the associator of the law $\mu$ and
$A_{\mu}^L (x_1\otimes x_2 \otimes x_3)=\mu(\mu(x_1, x_2),x_3)$ and $A_{\mu}^R =A_{\mu}^L -A_{\mu}$. As $A$ (resp. $B$)
is a $\mathcal{V}$-algebra, there exits $v \in \mathbb{K}$ (resp. $w \in \mathbb{K}$) such that $A_{\mu_A} \circ \Phi_v=0$
(resp. $A_{\mu_B} \circ \Phi_w=0)$. Let us suppose that $A \otimes B$ is a $\mathcal{V}-algebra$. There exists
$v' \in \mathbb{K}$ such that $A_{\mu_{A \otimes B}} \circ \Phi_{v'}=0).$ By taking $v'=\sum_{i=1}^6 \gamma_i \sigma_i$
the last condition can be written
$$
\sum_{i=1}^6 \gamma_i [A_{\mu_{A \otimes B}} \circ \Phi_{\sigma_i}]=0
$$
which can also be written
$$
\sum_{i=1}^6 \gamma_i [A_{\mu_{A}}^L \circ \Phi_{\sigma_i} \otimes A_{\mu_{B}}^L\circ \Phi_{\sigma_i}
-A_{\mu_{A}}^R\circ \Phi_{\sigma_i} \otimes A_{\mu_{B}}^R\circ \Phi_{\sigma_i}]=0.
$$
Let us denote
$e_i=A_{\mu_{A}}^L\circ \Phi_{\sigma_i},\widetilde{e_i}=A_{\mu_{A}}^R\circ \Phi_{\sigma_i},f_i=A_{\mu_{B}}^L\circ \Phi_{\sigma_i}$ and $\widetilde{f_i}=A_{\mu_{B}}^R\circ \Phi_{\sigma_i}.$ The vectors $e_i$ and $\widetilde{e_i}$ belongs to $Hom(A ^{\otimes ^3}, A)$ and the vectors $f_i$ and $\widetilde{f_i}$ to $Hom(B ^{\otimes ^3}, B)$.  The previous equation becomes:
$$(\star) \quad \sum_{i=1}^6 \gamma_i [e_i \otimes f_i - \widetilde{e_i} \otimes \widetilde{f_i}]=0.
$$
From the definition of the algebras $A$ and $B$ we have:
$$\sum_{i=1}^6 a_i (e_i -\widetilde{e_i})=0 \ {\rm and} \ \sum_{i=1}^6 b_i (f_i -\widetilde{f_i})=0$$
if $v=\sum_{i=1}^6 a_i \sigma_i$ and $w=\sum_{i=1}^6 b_i \sigma_i.$ Suppose that $dim F_w=k$ with $0\leq k \leq 6$. Then the rank of the vectors $\left\{ f_i,\widetilde{f_i} \right\}$ is equal to $(6+k).$ We can suppose that
$\left\{ f_1,...,f_6,\widetilde{f_1},...,\widetilde{f_k} \right\}$ are independent. This implies:
$$
\left\{
\begin{array}{l}
\widetilde{f}_{k+1}=\rho_1^{k+1}f_1+...+\rho_6^{k+1}f_6+\widetilde{\rho}_1^{k+1}\widetilde{f}_1+...+\widetilde{\rho}_k^{k+1}\widetilde{f}_k \\
: \\
: \\
\widetilde{f}_{6}=\rho_1^{6}f_1+...+\rho_6^{6}f_6+\widetilde{\rho}_1^{6}\widetilde{f}_1+...+\widetilde{\rho}_k^{6}
\widetilde{f}_k. \\
\end{array}
\right.
$$
The equation $(\star)$ is then written:
$$
\sum_{i=1}^6 e_i' \otimes f_i +\sum_{i=1}^k e_i" \otimes \widetilde{f}_k=0
$$
and the independence of the vectors $\left\{ f_1,...,f_6,\widetilde{f}_1,...,\widetilde{f}_k \right\}$ implies
$e'_1=...=e'_6=e"_1=...=e"_k=0.$ Then we deduce
$$
\left\{
\begin{array}{l}
\gamma_1 e_1- \gamma_{k+1}\rho_1^{k+1}\widetilde{e}_1-...-\gamma_{6}\rho_1^{6}\widetilde{e}_6=0 \\
: \\
: \\
\gamma_6 e_6- \gamma_{k+1}\rho_6^{k+1}\widetilde{e}_1-...-\gamma_{6}\rho_6^{6}\widetilde{e}_6=0 \\
\gamma_1 \widetilde{e_1}+ \gamma_{k+1}\widetilde{\rho}_1^{k+1}\widetilde{e}_{k+1}+...+\gamma_{6}\widetilde{\rho}_1^{6}\widetilde{e}_6=0 \\
: \\
: \\
\gamma_k \widetilde{e}_k+ \gamma_{k+1}\widetilde{\rho}_k^{k+1}\widetilde{e}_{k+1}+...+\gamma_{6}\widetilde{\rho}_k^{6}\widetilde{e}_6=0. \\
\end{array}
\right.
$$
Let us suppose that $k\neq 0.$ Then the conditions on the vector $v$ are of the type
$\sum_{i=1}^6 a_i(e_i -\widetilde{e}_i)=0$ and we have succesively $\gamma_1=...=\gamma_k=0$ and $\gamma_{k+1}\widetilde{\rho}_i^{k+1}=...=\gamma_{6}\widetilde{\rho}_i^{6}=0$ for $i=1,...,6$. If one of $\gamma_j$ for $j=k+1,...,6$ is not 0 then $\widetilde{\rho}_j^{i}=0$ for $i=1,...,6$. Thus
$$\widetilde{f}_j=\rho_1^{j}f_1+...+\rho_6^{j}f_6$$
and, as for every $i$ there exists $\sigma \in \Sigma _3$ such that $\widetilde{f}_i=\widetilde{f}_j \circ \Phi _{\sigma}$, we have also
$$\widetilde{f}_i=\rho_1^{i}f_1+...+\rho_6^{i}f_6$$
for every $i=1,..,6$. This implies $k=0$ which is impossible from the hypothesis.
Thus $k=0$ and the vectors $\left\{ f_i, \widetilde{f}_i \right\}_{i=1,...,6}$ is of rank 6.
The only possible relations are then $f_i=\widetilde{f}_i$ and $B$ is associative.
We deduce the associativity of $A.$

\medskip

\noindent{\bf Remark. } We have seen that a $W$-algebra is defined by the condition $A_{\mu}(x,x,x)=0$ for all $x$. Let ${\cal{A}}_1$ and ${\cal{A}}_2$ two $W$-algebras. For all $a \in {\cal{A}}_1$ and $b \in {\cal{A}}_2$ we have
$$
\begin{array}{l}
A_{\mu _1  \otimes \mu _2}(a\otimes b,a\otimes b,a\otimes b)= \\
 =\mu _1(\mu _1(a,a),a) \otimes \mu _2(\mu _2(b,b),b) -
\mu _1(a,\mu _1(a,a) \otimes \mu _2(b,\mu _2(b,b)) \\
=(\mu _1(\mu _1(a,a),a)-\mu _1(a,\mu _1(a,a))\otimes \mu _2(\mu _2(b,b),b) \\
=0.
\end{array}
$$
Then every indecomposable tensor in $A \otimes B$ generates an associative subalgebra. This doesnot imply, in general, that $A \otimes B$ is power-associative.

\section{A generalization : $\cal{V}$-$\cal{W}$-algebras}

The study of tensor products of $\mathcal{V}$-algebras has shown the necessity to write the associator as the difference
$$ A_{\mu}=A_{\mu}^L-A_{\mu}^R
$$
where $A_{\mu}^L(x_1,x_2,x_3)=\mu(\mu(x_1,x_2),x_3)$ and $A_{\mu}^R(x_1,x_2,x_3)=\mu(x_1,\mu(x_2,x_3)).$

Now, instead of considering action of $\Sigma_3$-permutation on the associator we can consider it independently on $A_{\mu}^L$ and $A_{\mu}^R$ which will induce different symmetries.

\subsection{Definition}

\begin{definition}
The catergory of $\cal{V}$-$\cal{W}$-algebras correspond to the $\mathbb{K}$-algebras $(\mathcal{A}, \mu)$ for which exist $v,w \in \mathbb{K}\,[\Sigma_3]$ such that
$$(\star)
\left\{
\begin{array}{l}
A_{\mu}^L \circ \Phi_v=0 \\
A_{\mu}^R \circ \Phi_w=0
\end{array}
\right.
\quad {\rm or} \quad (\star\star)
\left\{
\begin{array}{l}
A_{\mu}^L \circ \Phi_v-A_{\mu}^R \circ \Phi_w=0
\end{array}
\right.
$$
\end{definition}

We will write that $(\mathcal{A},\mu)$ is a $(v,w)$-algebra if $(\star)$ or $(\star \star)$ is satisfied and if for every $v',w'$
such that $v \in F_{v'}, \ w \in F_{w'}$ and $v' \notin F_v,\ w' \notin F_w$ we have $A_{\mu}^L \circ \Phi_{v'} \neq 0$
or $A_{\mu}^R \circ \Phi_{w'}\neq 0.$ It is clear that if $v=w$ or if $w \in F_v$ then a $(v,w)$-algebra is
a $v$-algebra. Now the problem is to look if for a $(v,w)$-algebra there exists $v'$ such that this algebra is a
$v'$-algebra. The most interesting example concerns the case of a $(v,w)$-algebra is a pre-Lie algebra or a $V$-algebra.

\subsection{Lie admissible $\cal{V}$-$\cal{W}$-algebras}

\begin{proposition}
Let $(\mathcal{A},\mu)$ be a $\cal{V}$-$\cal{W}$-algebras of type $(\star\star)$.
Then if $\cal A$ is a Lie-admissible algebra we have :

\noindent 1) $V \in F_V \cap F_W$

\noindent 2) One of the two following conditions is satisfied :

a) $V=\sum a_i\sigma_i(v)=\sum a_i\sigma_i(w)$ i.e $\exists \chi \in \mathbb{K}\,[\Sigma_3]$ such that
$\chi(v)=\chi(w)=V$,

b) there exits $\chi \in \mathbb{K}\,[\Sigma_3]$ such that $\chi (v-w)=0$ and $\chi (V) \neq 0.$

\end{proposition}

\noindent {\it Proof.} As the law $\mu$ satisfies
$$A_{\mu}^L \circ \Phi_v-A_{\mu}^R \circ \Phi_w=0
$$
and if such a law is Lie-admissible then we have that $V\in F_v$ and $V \in F_w.$
Let us first notice that to study $\cal{V}$-$\cal{W}$-algebras is reduced to study the representations of $\Sigma_3$
in the vector space of $\widetilde{E}$ generated by $(x_i.x_j).x_k$ and $x_i.(x_j.x_k)$ with
$i,j,k \in \left\{ 1,2,3 \right\}, \ i \neq j \neq k \neq i.$ This vector space is 12-dimensional and isomorphic to
$\mathbb{K}\,[\Sigma_3] \otimes \mathbb{K}\,[\Sigma_3].$ The invariant subvector spaces of $\widetilde{E}$ are also
generated by a vector $\widetilde{v}.$ When the algebra is a $\mathcal{V}$-algebra (considered as a
$\cal{V}$-$\cal{W}$-algebra) we restrict the study to the representations of $\Sigma_3$ in $\widetilde{E}$ which
let invariant the subspace $S=\left\{  (x_i.x_j).x_k-x_i.(x_j.x_k) \right\}.$ In this case we have a representation with
an interpretation as a representation of $\Sigma_3$ in $\mathbb{K}\,[\Sigma_3];$ this was the subject of the first part
of this work.

\noindent Let us come back to
$$A_{\mu}^L \circ \Phi_v-A_{\mu}^R \circ \Phi_w=0.$$
This relation can also be written
$$A_{\mu}^L \circ \Phi_v-A_{\mu}^R \circ \Phi_v+A_{\mu}^R \circ \Phi_v-A_{\mu}^R \circ \Phi_w=0
$$
or
$$A_{\mu} \circ \Phi_v-A_{\mu}^R \circ \Phi_u=0$$
with $u=w-v.$ As $V \in F_v$, if $\left\{ \sigma_1(v)=v, \sigma_2(v), ... , \sigma_k(v) \right\}$ ($k \leq 6$), is a basis of $F_v$, then $V=\sum a_i\sigma_i(v)$ and
$$A_{\mu} \circ \Phi_V = A_{\mu} \circ \Phi_{\sum a_i\sigma_i(v)}=\sum a_i A_{\mu} \circ \Phi_{\sigma_i(v)}.$$
From the definition of $\Phi_{\sigma}$ as linear transformation of ${\cal{A}} ^{\otimes ^3}$,$\Phi_{\sigma (v)}=\Phi _v \circ \Phi _{\sigma ^{-1}}$ which gives
$$A_{\mu} \circ \Phi_V=\sum a_i A_{\mu} \circ \Phi_v \circ \Phi _{\sigma _i ^{-1}}=A_{\mu} \circ \Phi_v \circ
\Phi _{\sum a_i \sigma _i ^{-1}}.$$
We deduce
$$A_{\mu} \circ \Phi_V -A_{\mu}^R \circ  \Phi_u \circ \Phi _{\sum a_i \sigma _i ^{-1}}=0.$$
Then
$$A_{\mu} \circ \Phi_V -A_{\mu}^R \circ  \Phi _{\sum a_i \sigma _i (u)}=0.$$
If the vector $u'=\sum a_i \sigma _i (u)$ is not reduced to the identity, $dim \ F_{u'} <6.$ If $F_{u'}=0$, then
$A_{\mu} \circ \Phi_V=0$ and the law $\mu$ is Lie-admissible. In this case $\sum a_i \sigma _i (u)=0$ that is
$\sum a_i \sigma _i (v)=\sum a_i \sigma _i (w)=V.$

\noindent If $dim \ F_{u'} \neq 0$, as $dim \ F_{u'} <6$, the dual space is not trivial. Thus there exists $\chi ' \in \Bbb{K} [\Sigma _3]$ such that $\chi' (u')=0$. The relation
$$A_{\mu} \circ \Phi_V-A_{\mu}^R \circ \Phi_{u'}=0$$
implies
$$A_{\mu} \circ \Phi_{\chi '(V)}-A_{\mu}^R \circ \Phi_{\chi '(u')}=0,$$
from which
$$A_{\mu} \circ \Phi_{\chi '(V)}=0.$$
But $\chi '(V)=aV.$ If $a \neq 0$ then $A_{\mu} \circ \Phi_V=0$ and the $\cal{V}-\cal{W}$-algebra is Lie-admissible. Let us put $\chi '= \sum c_i \sigma _i.$ Then
$$\chi '(u')=\sum c_i \sigma _i(u')  =
\sum c_ia_j \sigma _i(\sigma_j(u))=\sum c_ia_j (\sigma _j \circ \sigma _i)(u).$$ The vector $\chi=\sum c_ia_j\sigma _j \circ \sigma _i$ satisfies
$$\chi (u)=0.$$
$\clubsuit$

\begin{corollary}
If there exists $\chi ' \in \Bbb{K} [\Sigma _3]$ such that $\chi' (u')=0$ and $\chi '(V) \neq 0$ where $u'=\sum
_i \sigma _i (w-v)$ and $V=\sum a_i \sigma _i (v)$ then the $v-w$-algebra is Lie-admissible.
\end{corollary}

\medskip

\noindent {\it Examples}

1. Consider the $(v,w)$-algebra defined by
$$ x.(y.z)-x.(z.y)-(x.y).z+(y.x).z=0$$
Here $v=Id-\tau_{23}$ and $w=Id-\tau_{12}.$ We have
$$
\begin{array}{l}
V=v-\tau_{12}(v)-\tau_{13}(v)=\sum a_i \sigma_i(v) \\
V=w-\tau_{13}(w)-\tau_{23}(w)
\end{array}
$$
Then $V \in F_v \cap F_w.$ The vector $u=(v-w)$ is written
$$ w-v=(0,-1,0,1,0,0)
$$
We have
$$
\begin{array}{l}
\sum a_i\sigma_i(u)=u-\tau_{12}(u)-\tau_{13}(u)=(1,-1,0,1,0,-1)=u'
\end{array}
$$
Let $\chi'=Id +\tau_{12}+\tau_{13}.$ Then $\chi'(u')=0$ with $\chi=Id+c_1+c_2.$
This example plays a particular role in the study of Hopf operad. We will come back on this problem later.

2. Consider the $(v,w)$-algebra defined by
$$2(x.y).z-(y.x).z-(z.y).x-(x.z).y+(y.z).x-x.(y.z)-y.(z.x)-z.(x.y)
$$
Now $v=2Id-\tau_{12}-\tau_{13}-\tau_{23}+c_1$ and $w=Id+c_1+c_2.$
Here $dim F_v =5$ $dim F_w=2.$ The vector $u=w-v=(-1,1,1,1,0,1)$
We have $V=w-\tau_{12}(w).$ Let us take $u'=u-\tau_{12}(u)=(2,-2,0,-1,1,0).$
Then $dim F_{u'}=3$ and we have the relation $u' +\tau_{12}(u')=\tau_{23}(u')+c_1(u')=\tau_{13}(u')+c_2(u').$
This implies that every $\chi'$ such that $\chi'(u')=0$ has the following form $\chi'=(a_1,a_1,a_3,a_2,a_2,a_3).$
As $\chi'(v)=0$ the $(v,w)$-algebra is not Lie-admissible.
We can note that $u(V)$ is not zero.

\subsection{Monoidal category of $(\mathcal{V},\mathcal{W})$-algebras}

Contrary to the case of $\mathcal{V}$-algebras where the tensor product works only with two associative algebras,
we find by considering $(\mathcal{V},\mathcal{W})$-algebras, new examples of categories stable by tensor product.

\begin{proposition}
Let $\mathcal{C}$ be the category of $(\mathcal{V},\mathcal{W})$-algebras whose product satisfies:
$$\mu (\mu(x_1,x_2),x_3)-\mu(x_2,\mu(x_1,x_3))=0
$$
where $ x_1,x_2,x_3$ are elements of the algebra.
The category $\mathcal{C}$ is monoidal.
\end{proposition}

\noindent {\it Proof.} Let $A$ be an algebra. The algebra $A$ is a $(v,w)$-algebra by taking $v=Id$ and $w=c_1.$

Let $A$ and $B$ be two algebras of $\mathcal{C}$ and define the product of $A\otimes B$:
$$\forall \, a_i \in A, \forall \, b_j \in B, \quad (a_1\otimes b_1) \ ._{_{A\otimes B}} \ (a_2\otimes b_2)=(a_1 \ ._{_{A}} \ a_2) \otimes (b_1 \ ._{_{B}} \ b_2)
$$

The product defined in this way eqquipes the tensor product of these algebras with a structure of algebra of the same
type. In fact

$$
\begin{array}{ll}
[(a_1\otimes b_1)  .  (a_2\otimes b_2)] .  (a_3\otimes b_3)-
(a_2\otimes b_2)  .  [(a_1\otimes b_1)] .  (a_3\otimes b_3)] \\
=([a_1  .  a_2]  .  a_3) \otimes ([b_1  .  b_2]  .  b_3)-
(a_2  .  [a_1  . a_3]) \otimes (b_2  . [b_1  . b_3]) \\
=([a_1  .  a_2]  .  a_3) \otimes ([b_1  .  b_2]  .  b_3-b_2  . [b_1  . b_3]) \\
=0
\end{array}
$$

Then $A \otimes B \in \mathcal{C}. \qquad \clubsuit$

\section{Associated operads}

As we have natural right action of $\Sigma_3$ we can define trivially the binary quadratic operad associated
to each of these nonassociative type of algebras and their dual operads. We recall the basic definition of
a binary quadratic operad and its dual operad before giving them for the nonassociative algebras defined below.

\medskip

Let $\Bbb{K} [ \Sigma_n ]$ be the $\Bbb{K}$-group algebra of the symmetric group $\Sigma_n$.
An operad $\mathcal{P}$ is defined by a sequence of $\Bbb{K}$-vector spaces $\mathcal{P}%
(n)$, $n\geq 1$ such that $\mathcal{P}(n)$ is a right module over $\Bbb{K}\left[
\sum_{n}\right] $
with composition maps
\[
\circ _{i}:\mathcal{P}(n)\otimes \mathcal{P}(m)\rightarrow \mathcal{P}%
(n+m-1)\qquad i=1,...,n
\]
satisfying some ''associative'' properties, the May Axioms (\cite{M}, \cite{MSS}).

\noindent Any $\Bbb{K}\left[ \Sigma_{2}\right] $-module $E$ generates a free operad noted $\mathcal{F(}E)$
satisfying $\mathcal{F(}E)(1)=\Bbb{K}$, $\mathcal{F(}E)(2)=E$. In particular if $E=\Bbb{K}\left[ \Sigma_{2}\right]$, the free module $\mathcal{F(}E)(n)$ admits as a basis the
"parenthized products" of $n$ variables indexed by $\{1,2,...,n\}$. For instance a basis of
$\mathcal{F(}E)(2)$ is given by $(x_{1}.x_{2})$ and $(x_{2}.x_{1})$, and a basis of
$\mathcal{F(}E)(3)$ is given by
$$\left\{ \left( (x_{i}.x_{j}).x_{k}\right)
,\left( x_{i}.(x_{j}.x_{k}).\right) ,i\neq j\neq k\neq i,i,j,k\in \left\{
1,2,3\right\} \right\} .$$
Let $E$ be a $\Bbb{K}\left[ \Sigma_{2}\right]$-module and $R$ a $\Bbb{K}\left[ \Sigma_{3}\right]$-submodule of
$\mathcal{F(}E)(3)$. We denote $\mathcal{R}$ the ideal generated by $R$, that is the intersection of all the ideals
$\mathcal{I}$ of $\mathcal{F}(E)$ such that $\mathcal{I}(1)=0,\mathcal{I}(2)=0$ and
$\mathcal{I}(3)=R$.

We call binary quadratic operad generated by $E$ and $R$ the operad $\mathcal{P}(\Bbb{K},E,R)$, also
denoted  $\mathcal{F(}E)/\mathcal{R}$ and defined by
\[
\mathcal{P}(\Bbb{K},E,R)(n)=
(\mathcal{F(}E)/\mathcal{R)}(n)=\frac{\mathcal{F(}E)(n)}{\mathcal{R(}n%
\mathcal{)}}
\]

\noindent Thus an operad $\mathcal{P}$ is binary quadratic operad if and only if there exists a $%
\Bbb{K}\left[ \Sigma_{2}\right] $-module $E$ and a $\Bbb{K}\left[
\Sigma_{3}\right] $-submodule $R$ of $\mathcal{F(}E)(3)$ such that $%
\mathcal{P\simeq F(}E)/\mathcal{R}$.

\smallskip

\noindent {\bf Examples}.

1. The associative operad $\mathcal{A}ss,$ the Lie operad $\mathcal{L}%
ie$ \cite{G.K}.

2. Let $\mathcal{F(}E)$ be the free operad generated by $E=\Bbb{K}\left[ \Sigma_{2}\right]$.
Consider the $\Bbb{K}\left[ \Sigma_{3}\right] $-submodule $R$ generated by the
vector
\begin{eqnarray*}
u
&=&x_{1}.(x_{2}.x_{3})+x_{2}.(x_{3}.x_{1})+x_{3}.(x_{1}.x_{2})-x_{2}.(x_{1}.x_{3})-x_{3}.(x_{2}.x_{1}) \\
&&-x_{1}.(x_{3}.x_{2})
-(x_{1}.x_{2}).x_{3}-(x_{2}.x_{3}).x_{1}-(x_{3}.x_{1}).x_{2}+(x_{2}.x_{1}).x_3 \\
&& +(x_{3}.x_{2}).x_{1}+(x_{1}.x_{3}).x_{2}
\end{eqnarray*}

\noindent From now we will take
$$E=\Bbb{K}\left[ \Sigma_{2}\right].$$

\medskip

\noindent The Lie-Admissible operad, denoted $\mathcal{L}ieAdm$ is the binary quadratic operad
defined by
\[
\mathcal{L}ieAdm=\mathcal{F(}E)/\mathcal{R}.
\]

\medskip

For each of the above types of Lie-admissible and power-associative algebras there exists the
corresponding operad. The way we define these algebras gives directely the vectors generating,
as $\Bbb{K}\left[ \Sigma_{3}\right] $-submodule, the module $R$ and thus the binary quadratic algebra
associated to these algebras. Considering a $v$-Lie-admissible algebra, the corresponding module $F_v$ which determine
the full class of this $v$-algebra, has a binary quadratic operad whose corresponding module $R$ is $F_v$. From the
previous theorem, we know for each class of $v$-algebra the module of the associated operad. Then
 the relations of definition for all the operads of $v$-Lie-admissible algebras are given by this theorem.

\subsection{The dual operads associated to the Lie-admissible case}
Let us consider the binary quadratic operad
$\mathcal{P}(\Bbb{K},E,R)$. Then the dual binary quadratic
operad is defined by
$$\mathcal{P}^{!}=\mathcal{P}(\Bbb{K},E^{\vee },R^{\perp })$$
where $E^{\vee }$ is the dual of $E$ tensored by the signum representation of $\Sigma_n$ and $R^{\perp }$ the orthogonal
complement to $R$ in $\mathcal{F(}E^{\vee })(3)=\mathcal{F(}E)(3)^{\vee }$.

\begin{proposition}
The dual operads of algebras which are Lie-admissible are quadratic operads
whose corresponding algebras are associative algebras
satisfying respectively :

\medskip

\noindent - for type $(I)$ : $x_1.x_2.x_3=x_{\sigma(1)}.x_{\sigma(2)}.x_{\sigma(3)}$ for all $\sigma \in \sigma_3.$

\medskip

\noindent - for type $(II)$ : $x_1.x_2.x_3=x_2.x_3.x_1=x_3.x_1.x_2$

\medskip

\noindent - for type $(III_1)$ : $x_1.x_2.x_3+ tx_2.x_1.x_3-x_1.x_3.x_2-tx_3.x_1.x_2=0$
with $t \neq 1$

\medskip

\noindent - for type $(III_2)$ : $x_1.x_2.x_3= x_3.x_2.x_1$

\medskip

\noindent - for type $(III_3)$ : $x_1.x_2.x_3-x_2.x_1.x_3-2x_1.x_3.x_2+2x_2.x_3.x_1=0$

\medskip

\noindent -for type $(IV_1)$ : $t \neq 1$
$$
\begin{array}{ll}
(t-1) x_1.x_2.x_3-(t-1) x_2.x_1.x_3-(t+2) x_3.x_2.x_1+(1+2t) x_1.x_3.x_2 \\
-(1+2t)x_2.x_3.x_1+ (t+2) x_3.x_1.x_2=0
\end{array}
$$

\noindent -for type $(IV_2)$ : $x_1.x_2.x_3+x_2.x_1.x_3-x_3.x_2.x_1-x_3.x_1.x_2=0 $

\noindent -for type $(IV_3)$ : $x_1.x_2.x_3+x_2.x_1.x_3-x_1.x_3.x_2-x_3.x_1.x_2=0 $

\noindent -for type $(V)$ : $x_1.x_2.x_3-x_2.x_1.x_3-x_3.x_2.x_1-x_1.x_3.x_2+x_2.x_3.x_1+x_3.x_1.x_2=0$

\end{proposition}

\noindent {\it Proof.}  Let us consider the scalar product on $\mathcal{F(}E)(3)$ defined by
\begin{eqnarray*}
&<i(jk),i(jk)>=sgn(
\begin{tabular}{lll}
$1$ & $2$ & $3$ \\
$i$ & $j$ & $k$%
\end{tabular}
) \\
&<(ij)k,(ij)k>=-sgn(
\begin{tabular}{lll}
$1$ & $2$ & $3$ \\
$i$ & $j$ & $k$%
\end{tabular}
)
\end{eqnarray*}
Let $R$ be the $\Bbb{K}\left[ \Sigma_{3}\right] $-submodule which determines the
operad of algebras of type $(I)$. This corresponds to the Lie-admissible algebras or
$G_6$-algebras in \cite{G.R} hence the result.

\smallskip

If $R$ determines the
operad of algebras of type $(II)$ we get the $G_5$-algebras of \cite{G.R}.

\smallskip

Let $R$ be associated to one of the $(III_i)$ case.
  The annihilator $R^{\perp }$ respect to this scalar
product is of dimension $9$. Let $R^{\prime }$ be the $\Bbb{K}\left[
\Sigma_{3}\right] $-submodule of $\mathcal{F(}E)(3)$ generated by the
vectors $(x_{1}x_{2})x_{3}-x_{1}(x_{2}x_{3})$ and

-for i=1  $(x_{1}.x_{2}).x_{3}+t(x_2.x_1).x_{3}-(x_{1}.x_{3}).x_{2}-t(x_2.x_3).x_{1}$ with $t \neq 1;$

-for i=2  $(x_{1}x_{2})x_{3}-(x_{3}x_2)x_1;$

-for i=3 $(x_1.x_2).x_3+(x_3.x_2).x_1-(x_2.x_3).x_1+(x_2.x_1).x_3.$

Then $\dim \ R^{\prime }=9$ and $<u,v>=0$ for all $v\in R^{\prime }$ where $u$
is the vector which generates $R.$ This implies $R^{\prime }\simeq R^{\perp }
$ and ($\mathcal{F(}E)/\mathcal{R})^{!}$ is by definition the quadratic operad
$\mathcal{F(}E)/\mathcal{R}^{\perp }.$

\smallskip

We demonstrate in the same way the case of $(IV_i)$ with
$R^{\perp }$ is the
$\Bbb{K}\left[ \sum_{3}\right] $-sub-module of
$\mathcal{F(}E)(3)$ generated by the vectors
$(x_{1}x_{2})x_{3}-x_{1}(x_{2}x_{3})$ and

-for i=1 $(t-1) (x_1.x_2).x_3-(t-1) (x_2.x_1).x_3-(t+2) (x_3.x_2).x_1+$

$(1+2t) (x_1.x_3).x_2-(1+2t)(x_2.x_3).x_1+ (t+2) (x_3.x_1).x_2$ with $t \neq 1$.

-for i=2  $(x_{1}.x_{2}).x_{3}+(x_2.x_1).x_{3}-(x_{3}.x_{2}).x_{1}-(x_3.x_1).x_{2};$

-for i=3  $(x_{1}x_{2})x_{3}+(x_{2}x_1)x_3-(x_{1}x_{3})x_{2}-(x_{3}x_1)x_2.$

For the case $(V)$, the sub-module $R$ is of dimension 7 and generated as $\Sigma_3$-submodule
by $(x_{1}x_{2})x_{3}-x_{1}(x_{2}x_{3})$ and
$(x_{1}x_{2})x_{3}-(x_{2}x_{1})x_{3}-(x_{3}x_{2})x_{1}-(x_{1}x_{3})x_{2}+(x_{2}x_{3})x_{1}+(x_{3}x_{1})x_{2}.$

\subsection{The dual operads associated to the power-associative case}

The determination of the dual operads in the power-associative case is similar to the Lie-admissible one. When $F_v$ is 2 or 4-dimensional, as the $v$-algebras are both power-associative and Lie-admissible and
the dual operads are described in the previous proposition.

\begin{proposition}
The dual operads of power-associative algebras are quadratic operads
whose corresponding algebras are associative algebras
satisfying respectively :

\medskip

\noindent - for type $(I')$ : $x_1.x_2.x_3=-x_2.x_1.x_3=-x_1.x_3.x_2 $

\medskip

\noindent - for type $(III'_{\ 1})$ :

$-2x_1.x_2.x_3-(2+t)x_3.x_2.x_1+(t-1)x_1.x_3.x_2-(1+t)x_2.x_3.x_1+tx_3.x_1.x_2=0 $

\medskip

\noindent - for type $(III'_{\ 2})$ : $x_1.x_2.x_3+x_2.x_1.x_3+x_3.x_2.x_1+x_2.x_3.x_1=0 $

\medskip

\noindent - for type $(V')$ : $x_1.x_2.x_3+x_2.x_1.x_3+x_3.x_2.x_1+x_1.x_3.x_2+x_2.x_3.x_1+x_3.x_1.x_2=0$

\medskip

When $F_v$ is 2 or 4-dimensional, as the $v$-algebras are both power-associative and Lie-admissible, we have define
the dual operad in the previous proposition.
\end{proposition}

The proof is similar for the Lie-admissible and power-associative cases.

\subsection{An interpretation of algebras on the dual operad}

If we look the way to define the dual operad, we see that
$\mathcal{P}^{!}=\frac{\mathcal{F}(E)}{\mathcal{R}^{\perp}}$ with $\mathcal{R}^{\perp}=<(x_1.x_2).x_3-x_1.(x_2.x_3),u_1,...,u_p>.$
So $\mathcal{R}^{\perp}$ can be decomposed in two subspaces $R_{ass}$ which gives the associativity and a second one $U$ which gives the supplementary rules that an algebra on this operad has to satisfy.

\smallskip

\noindent {\bf Example.} In the case of Vinberg algebras, we have
$$
R^{\perp}=R_{ass} \oplus <(x_1.x_2).x_3-(x_2.x_1).x_3>.
$$
We can make such a decomposition so that every vector of $U$ is left parenthized and that $U$ is a
$\mathbb{K}\,[\Sigma_3]$-module. Let us consider the $\mathbb{K}\,[\Sigma_3]$-module
$\widetilde{E}=\left\{ x_i.x_j.x_k, i \neq j \neq k \neq i \right\}$ of non parenthized 3-products, $dim \widetilde{E}=6.$
Then $U$ is isomorphic to a submodule $G$ of $\widetilde{E}.$ There is a natural right action of
$\mathbb{K}\,[\Sigma_3]$ on $\widetilde{E}:$

$$
\begin{array}{ccc}
\widetilde{E} \times \mathbb{K}\,[\Sigma_3] & \rightarrow & \widetilde{E} \\
(x_1.x_2.x_3,\sigma ) & \rightarrow & x_{\sigma^{-1}(1)}.x_{\sigma^{-1}(2)}.x_{\sigma^{-1}(3)}.
\end{array}
$$

If we consider $A$ a $v$-algebra, the binary quadratic operad $\mathcal{P}$ associated to this algebra is such that
$R=F_v.$ We can naturally associated a sub-module $H$ of $\widetilde{E}$ to $F_v:$
$$
\begin{array}{ccc}
 E &  \rightarrow & \widetilde{E} \\
v=  A_{\mu}\circ \sigma{x_1 \otimes x_2 \otimes x_3} & \rightarrow & \widetilde{v}= x_{\sigma^{-1}(1)}.x_{\sigma^{-1}(2)}.x_{\sigma^{-1}(3)}
\end{array}
$$
Then in $\widetilde E$, $dim H=dim F_v$ and $dim G=codim H=6- dim H.$

We can also define a notion of $\widetilde{v}$-algebra in $\widetilde{E}$ that is $(A,\mu)$ is a
$\widetilde{V}$-algebra if and only if $\Phi_{\widetilde{v}}=0.$ We have similar results for the orbit of
$\widetilde{v}$ under the previous action and decomposition of $\widetilde{\mathcal{O}}(\widetilde{v})$
in direct sum of irreducible invariant subspaces.

Let $V$ and $W$ be the vectors of $\mathbb{K}\,[\Sigma_3]$ corresponding to Lie-admissible and power-associative.
As $<V,W>\neq 0,$ we have directely that if $V$ and $W$ are both in $F_v,$ neither $\widetilde{V}$ nor $\widetilde{W}$
in $R^{\perp}.$ Then $G$ is a direct sum of 2-dimensional irreducible invariant submodule of $E.$

\begin{proposition}
Let $A$ be a $v$-algebra over an operad $\mathcal{P}=\frac{\mathcal{F}(E)}{\mathcal{R}}.$

$\bullet$ If $F_v$ contains $V$ and $W$ then $R^{\perp}$ defining the dual operad is
$$R^{\perp}=R_{ass} \oplus_i F_{w_i}$$
where $F_{w_i}$ are 2-dimensional irreducible invariant submodules.

\smallskip

$\bullet$ If $F_v$ contains $V$ (resp. $W$) but not $W$ (resp. $V$) then
$$R^{\perp}=R_{ass} \oplus V \oplus_i F_{w_i}  \ ( {\rm resp.} \ R_{ass} \oplus W \oplus_i F_{w_i} )$$
\end{proposition}

\noindent {\it Example.} If we consider a Lie-algebra of type $(IV_1)$ defined by
$v=(t-1) Id-(t-1) \tau_{12}-(t+2) \tau{13}+(1+2t) \tau{23}-(1+2t)c_1+ (t+2) c_2, \ t \neq 1$, as $V$ and $W$ are in $F_v$,
we can verify that $R^{\perp}=U \oplus F_{w_1}$ where $F_{w_1}$ is a 2-dimansional irreducible space.

\section{On the tensor product of algebras on operads}

We have seen in \cite{G.R} that $G_i -ASS$ the category $G_i$-associative algebras on $\Bbb{K}$ for $i=2,...,6$ is not associated to a Hopf operad which implies that two algebras in the same category of $G_i$-algebras ($i \geq 6$)
do not give in general a $G_i$-algebra. But the functor $A \otimes -$ determines a one-to-one correspondence
between the algebras over the operad $G_{i}-\mathcal{A}ss$ and algebras on the dual operad $G_{i}-\mathcal{A}ss^{! }$.

\begin{theorem}[G.R]
Let $A$ be a $G_i$-associative algebra. Then  $A \otimes -$ is a covariant functor
$$ A \otimes - : (G_i -ASS)^{!} \longrightarrow G_i -ASS $$
where $ (G_i -ASS)^{!}$  is the category of algebras corresponding to the algebras over the dual operad $G_{i}-\mathcal{A}ss^{!}$.
\end{theorem}

Moreover, we can see that
\begin{theorem}
If $A$ is a $V$-algebra and $B$ a $W^!$-algebra (that is an algebra on the dual operad of the operad of W-algebras) then $A \otimes B$ can be provided with a product of $W$-algebra;
if $A$ is a $W$-algebra and $B$ a $V^!$-algebra then $A \otimes B$ is a $W$-algebra;
if $A$ is a $W$-algebra and $B$ a $W^!$-algebra then $A \otimes B$ is a $V$-algebra.
\end{theorem}

This shows a link between $1$-dimensional $\mathcal{V}$ and $\mathcal{W}$-algebras and their "dual" algebras.
In the case of $3$-dimensional $\mathcal{V}$ and $\mathcal{W}$-algebras, we have the same type of results that is if we call $\bar{G}_i-Ass$ and $\bar{G}_i-\mathcal{A}ss$ the corresponding category and operad associated to the ${G}_i-Ass$ and $G_i-\mathcal{A}ss$ in the power-associative case (that is $\bar{G}_1-Ass$ is the category of $W$-algebras, $\bar{G}_2-Ass$ is the category of $v$-algebras with $v=Id+\tau_{12}$, ...) we have that
\begin{theorem}
If $A$ is a $G_i-Ass$-algebra and $B$ a $\bar{G}_i^!-Ass$-algebra then $A \otimes B$ can be provided with a product of $G_i$-algebra;
If $A$ is a $\bar{G}_i-Ass$-algebra and $B$ a $G_i^!-Ass$-algebra then $A \otimes B$ is a $\bar{G}_i-Ass$-algebra;
If $A$ is a $\bar{G}_i$-algebra and $B$ a $\bar{G}_i^!-Ass$-algebra then $A \otimes B$ is a $\bar{G}_i-Ass$-algebra;
\end{theorem}

The description of operads in terms of $\Sigma_3$-invariant subspaces permits to generalize the case of $G_i$-algebras.
What about the tensorial construction for the new Lie-admissible and power-associative algebras?

In the previous examples, what enables the construction is that we can factorize the fist part to make appear the fact that the algebra $A \otimes B$ is a $\mathcal{V}$ or $\mathcal{W}$-algebra.
In fact, we have that $R{\perp}=R_{ass} \oplus U$ with $U$ isomorphic to $<Id+\varepsilon_i \sigma_i>$
with
$\varepsilon_i \in \left\{ -1,1 \right\}$ and $\sigma_i \in \sigma_3$, $\sigma_i \neq Id$.

The possibilities for $U$ are $<Id-\tau_{12}>,<Id-\tau_{13}>,<Id-\tau_{23}>$ which are of dimension 3 and correspond to the Vinberg, pre-Lie and $G_5$ case, $<Id+\tau_{12}>,<Id+\tau_{13}>,<Id+\tau_{23}>$  of dimension 3 and which are the corresponding power-associative cases, $<Id-c_1>$ which is of dimension 4 and correspond to the dual operad of $G_5(=\bar{G}_5)$, $<Id-\tau_{12}, Id-c_1>, <Id+\tau_{12}, Id-c_1>$ which are of dimension 5 and correspond respectively to the $G_6$ and $\bar{G}_6$ cases,  $<Id+c_1>$ which is of dimension 6 and corresponds to the associative case
($<Id+c_1>=<Id>$).

\end{document}